\documentclass[11pt,letterpaper]{article}
\usepackage{amsfonts, amsmath, amssymb, amscd, amsthm, color, graphicx, mathrsfs, wasysym, setspace, mdwlist}
\usepackage{hyperref}
\newcommand{\G}{\Gamma (G, X\sqcup H)}

\newcommand{\h}{\hookrightarrow_h}

\newcommand{\e}{\varepsilon}

\renewcommand{\d}{{\rm d}}

\renewcommand{\ll }{\langle\hspace{-.7mm}\langle }
\newcommand{\rr }{\rangle\hspace{-.7mm}\rangle }

\newtheorem{thm}{Theorem}[section]
\newtheorem{cor}[thm]{Corollary}

\newtheorem{q}[thm]{Question}

\theoremstyle{definition}
\newtheorem{defn}[thm]{Definition}
\theoremstyle{remark}

\newtheorem{ex}[thm]{Example}

\newfont{\eufm}{eufm10}

\begin{document}

\title{Groups acting acylindrically on hyperbolic spaces}

\author{D.V. Osin\thanks{The author was supported by the NSF grant DMS-1612473.}}

\date{}

\maketitle

\begin{abstract}
The goal of this article is to survey some recent developments in the study of groups acting on hyperbolic spaces. We focus on the class of \emph{acylindrically hyperbolic groups} and their \emph{hyperbolically embedded subgroups}. This class is broad enough to include many examples of interest, yet a significant part of the theory of hyperbolic and relatively hyperbolic groups can be generalized in this context. 
\end{abstract}

%\hspace*{3mm}\textbf{2010 MSC:} 20F65, 20F67, 20F69.

\tableofcontents

%%%%%%%%%%%%%%%%%%%%%%%%%%%%%%%%%%%%%%%%%%%%%%%%%%%%%%%%%%%%%%%%%%%%%%%

\section{Introduction}

%%%%%%%%%%%%%%%%%%%%%%%%%%%%%%%%%%%%%%%%%%%%%%%%%%%%%%%%%%%%%%%%%%%%%%%

Suppose that a group $G$ acts by isometries on a metric space $S$. If the action is sufficiently ``nice", many properties of $G$ can be revealed by studying the geometric structure of $G$-orbits in $S$. This approach works especially well if $S$ satisfies certain negative curvature condition. 

Systematic research in this direction began in late 1980s when Gromov \cite{Gro} introduced the notion of an abstract hyperbolic metric space.
Groups acting properly and cocompactly on hyperbolic spaces are called \emph{word hyperbolic}. More generally, replacing properness with its relative analogue modulo a fixed collection of subgroups leads to the notion of a \emph{relatively hyperbolic group}. The study of hyperbolic and relatively hyperbolic groups was initiated by \cite{Gro} and since then it has been one of the most active areas of research
in geometric group theory.

A further generalization, the class of \emph{acylindrically hyperbolic groups}, was suggested in \cite{Osi16a} and received considerable attention in the past few years. It includes many examples of interest: non-elementary hyperbolic and relatively hyperbolic groups, all but finitely many mapping class groups of punctured closed surfaces, $Out(F_n)$ for $n\ge 2$, most $3$-manifold groups, groups of deficiency at least $2$, and the Cremona group of birational transformations of the complex projective plane, just to name a few. On the other hand, the property of being acylindrically hyperbolic is strong enough to allow one to apply powerful geometric techniques.

A significant part of the theory of relatively hyperbolic groups can be generalized to acylindrically hyperbolic groups using the notion of a \emph{hyperbolically embedded collection of subgroups} introduced in \cite{DGO}. In particular, this notion provides a suitable framework for developing a group theoretic version of Thurston's theory of hyperbolic Dehn filling in $3$-manifolds. Group theoretic Dehn filling was originally studied in the context of relatively hyperbolic groups in \cite{GM,Osi07}. Recently it was used to obtain several deep results (most notably, it was employed in Agol's proof of the virtual Haken conjecture \cite{A}). Yet another powerful tool is \emph{small cancellation theory}, which can be used to prove various embedding theorems and to construct groups with unusual properties \cite{Hull,Osi10}.

The main purpose of this paper is to survey the recent progress in the study of acylindrically hyperbolic groups and their hyperbolically embedded subgroups. In the next section we briefly discuss equivalent definitions, main examples, and basic properties of acylindrically hyperbolic groups. Hyperbolically embedded subgroups are discussed in Section 3. Section 4 is devoted to group theoretic Dehn filling. An informal discussion of small cancellation theory and a survey of some application is given in Section 5.

%%%%%%%%%%%%%%%%%%%%%%%%%%%%%%%%%%%%%%%%%%%%%%%%%%%%%%%%%%%%%%%%%%%%%%%

\section{Acylindrically hyperbolic groups}\label{2}

%%%%%%%%%%%%%%%%%%%%%%%%%%%%%%%%%%%%%%%%%%%%%%%%%%%%%%%%%%%%%%%%%%%%%%%

\paragraph{2.1. Hyperbolic spaces and group actions.}
We begin by recalling basic definitions and general results about groups acting on hyperbolic spaces. Our main reference is the Gromov's paper \cite{Gro}; additional details can be found in \cite{BH} and \cite{GH}. All group actions on metric spaces discussed in this paper are assumed to be isometric by default.

\begin{defn}
A metric space $S$ is \emph{hyperbolic} if it is geodesic and there exists $\delta \ge 0$ such that for any geodesic triangle $\Delta $ in $S$, every side of $\Delta $ is contained in the union of the $\delta$-neighborhoods of the other two sides.
\end{defn}

\begin{ex}
Every bounded space $S$ is hyperbolic with $\delta ={\rm diam} (S)$. Every tree is hyperbolic with $\delta=0$. $\mathbb H^n$ is hyperbolic for every $n\in \mathbb N$. On the other hand,
 $\mathbb R^n$ is not hyperbolic for $n\ge 2$.
\end{ex}

Given a hyperbolic space $S$, we denote by $\partial S$ its \emph{Gromov boundary}. We do not assume that the space is proper and therefore the boundary is defined as the set of equivalence classes of sequences of points convergent at infinity; for details we refer to \cite[Section 1.8]{Gro}. The union $\widehat S= S\cup \partial S$ is a completely metrizable Hausdorff topological space containing $S$ as a dense subset.

\begin{ex}
The Gromov boundary of a bounded space is empty. $\partial (\mathbb H^n) =\mathbb S^{n-1}$. The boundary of an $n$-regular tree is the Cantor set if $n\ge 3$ and consists of two points if $n=2$.
\end{ex}

Let $G$ be a group acting (isometrically) on a hyperbolic space $S$. This action extends to an action on $\widehat S$ by homeomorphisms. We denote by $\Lambda (G)$ the \emph{limit set} of $G$, that is, the set of accumulation points of a $G$-orbit on $\partial S$. Thus $$\Lambda (G)= \overline{Gs}\cap \partial S,$$ where $s\in S$ and $\overline{Gs}$ is the closure of the corresponding orbit.  In fact, this definition is independent of the choice of $s\in S$. Given an element $g\in G$, we denote $\Lambda (\langle g\rangle)$ simply by $\Lambda (g)$ and call it the \emph{limit set} of $g$.

Similarly to the classification of elements of $PSL(2,\mathbb R)= \mathrm{Isom}\,\mathbb H^2$, we have the following classification of isometries of abstract hyperbolic spaces.

\begin{defn}
An element $g\in G$ is called \emph{elliptic} if $\Lambda (g)=\emptyset$ (equivalently, all orbits of $\langle g\rangle $ are bounded), \emph{parabolic}  if $|\Lambda (g)|=1$; and \emph{loxodromic} if $|\Lambda (g)|=2$.
Equivalently, an element $g\in G$ is loxodromic if the map $\mathbb Z\to S$ defined by $n\mapsto g^ns$ is a quasi-isometric embedding for every $s\in S$; in turn, this is equivalent to the existence of $c>0$ such that $\d_S(s, g^ns)\ge c|n|$ for all $n\in \mathbb Z$. Two loxodromic elements $g,h\in G$ are called \emph{independent} if $\Lambda(g)\cap \Lambda (h)=\emptyset$.
\end{defn}

We recall the standard classification of groups acting on hyperbolic spaces, which goes back to Gromov \cite[Section 8.2]{Gro}.

\begin{thm}[Gromov]\label{class1} For every group $G$ acting on a hyperbolic space $S$, exactly one of the following conditions holds.
\begin{enumerate}
\item[1)] $|\Lambda (G)|=0$. Equivalently,  $G$ has bounded orbits. In this case the action of $G$ is called \emph{elliptic}.

\item[2)] $|\Lambda (G)|=1$. Equivalently, $G$ has unbounded orbits and contains no loxodromic elements. In this case the action of $G$ is called \emph{parabolic}.

\item[3)] $|\Lambda (G)|=2$. Equivalently, $G$ contains loxodromic elements and any two loxodromic elements have the same limit points. In this case the action of $G$ is called \emph{lineal}.

\item[4)] $|\Lambda (G)|=\infty$. Then $G$ always contains loxodromic elements. In turn, this case breaks into two subcases.
\begin{enumerate}
\item[a)] $G$ fixes a point $\xi \in \partial S$. In this case $\xi$ is the common limit point of all loxodromic elements of $G$. Such an action is called \emph{quasi-parabolic}.
\item[b)] $G$ has no fixed points on $\partial S$. Equivalently, $G$ contains independent loxodromic elements. In this case the action is said to be of \emph{general type}.
\end{enumerate}
\end{enumerate}
\end{thm}

\begin{defn}
The action of $G$ is called \emph{elementary}  in cases 1)--3) and \emph{non-elementary} in case 4).
\end{defn}

An action of a group $G$ on a metric space $S$ is called (metrically) \emph{proper} if the set $\{ g\in G\mid \d_S (s, gs)\le r\}$ is finite for all $s\in S$ and $r\in \mathbb R_+$. Further, the action of $G$ is \emph{cobounded}  if there exists a bounded subset $B\subseteq S$ such that $S=\bigcup_{g\in G} gB$. Finally, the action is \emph{geometric} if it is proper and cobounded. (We work in the category of metric spaces here, so compactness gets replaced by boundedness.)

For geometric actions, we have the following, see \cite{Gro}.

\begin{thm}[Gromov]\label{class2}
Let $G$ be a group acting geometrically on a hyperbolic space. Then exactly one of the following three conditions hold.
\begin{enumerate}
\item[(a)] $G$ acts elliptically. In this case $G$ is finite.
\item[(b)] $G$ acts lineally. In this case $G$ is virtually cyclic.
\item[(c)] The action of $G$ is of general type.
\end{enumerate}
\end{thm}

To every group $G$ generated by a set $X$ one can associate a natural metric space, namely the Cayley graph $\Gamma (G,X)$, on which $G$ acts geometrically. The vertex set of $\Gamma (G,X)$ is $G$ itself and two elements $g,h$ are connected by an edge if $g=hx$ for some $x\in X^{\pm 1}$. This graph is endowed with the \emph{combinatorial metric} induced by identification of edges with $[0,1]$.

\begin{defn}
A group $G$ is \emph{hyperbolic} if it admits a geometric action on a hyperbolic space.
\end{defn}

Equivalently, a group $G$ generated by a finite set $X$ is hyperbolic if the Cayley graph $\Gamma (G,X)$ is a hyperbolic metric space. The equivalence of these two definitions follows from the well-known Svarc-Milnor Lemma and quasi-isometry invariance of hyperbolicity of geodesic spaces, see \cite{BH, Gro} for details.

\paragraph{2.2. Equivalent definitions of acylindrical hyperbolicity.}
Recall that the action of a group $G$ on a metric space $S$ is {\it acylindrical} if for every $\e>0$ there exist $R,N>0$ such that for every two points $x,y$ with $\d (x,y)\ge R$, there are at most $N$ elements $g\in G$ satisfying
$$
\d(x,gx)\le \e \;\;\; {\rm and}\;\;\; \d(y,gy) \le \e.
$$
The notion of acylindricity goes back to Sela's paper \cite{Sel}, where it was considered for groups acting on trees. In the context of general metric spaces, the above definition is due to Bowditch \cite{Bow}. Informally, one can think of this condition as a kind of properness of the action on $S\times S$ minus a ``thick diagonal".

\begin{ex}\label{aact}
\begin{enumerate}
\item[(a)] If $S$ is a bounded space, then every action $G\curvearrowright S$ is acylindrical. Indeed it suffices to take $R>{\rm diam}(S)$.
\item[(b)] It is easy to see that every geometric action is acylindrical. On the other hand, proper actions need not be acylindrical in general.
\end{enumerate}
\end{ex}

We begin with a classification of groups acting acylindrically on hyperbolic spaces. The following theorem is proved in \cite{Osi16a} and should be compared to Theorems \ref{class1} and \ref{class2}

\begin{thm}\label{class3}
Let $G$ be a group acting acylindrically on a hyperbolic space. Then exactly one of the following three conditions holds.
\begin{enumerate}
\item[(a)] $G$ acts elliptically, i.e., $G$ has bounded orbits.
\item[(b)] $G$ acts lineally. In this case $G$ is virtually cyclic.
\item[(c)] The action of $G$ is of general type.
\end{enumerate}
\end{thm}

Compared to the general classification of groups acting on hyperbolic spaces, Theorem \ref{class3} rules out parabolic and quasi-parabolic actions and characterizes groups acting lineally. On the other hand, compared to Theorem \ref{class2}, finiteness of elliptic groups is lacking. This part of Theorem \ref{class3} cannot be improved, see Example \ref{aact} (a).

Applying the theorem to cyclic groups, we obtain the following result first proved by Bowditch \cite{Bow}.

\begin{cor}
Every element of a group acting acylindrically on a hyperbolic space is either elliptic or loxodromic.
\end{cor}

\begin{defn}\label{defah}
We call a group $G$ \emph{acylindrically hyperbolic} if it admits a non-elementary acylindrical action on a hyperbolic space. By Theorem \ref{class3}, this is equivalent to the requirement that $G$ is not virtually cyclic and admits an acylindrical action on a hyperbolic space with unbounded orbits.
\end{defn}

Unfortunately, Definition \ref{defah} is hard to verify in practice. Instead, one often first proves that the group satisfies a seemingly weaker condition, which turns out to be equivalent to acylindrical hyperbolicity. To formulate  this condition we need a notion introduced by Bestvina and Fujiwara in \cite{BF}.

\begin{defn}\label{WPD}
Let $G$ be a group acting on a hyperbolic space $S$, $g$ an element of $G$.  One says that $g$ satisfies the {\it weak proper discontinuity} condition (or $g$ is a {\it WPD element}) if for every $\e >0$ and every $s\in S$, there exists $M\in \mathbb N$ such that
\begin{equation}\label{eq: wpd}
\left| \{ a\in G \mid \d _S(s, as)<\e, \;   \d (g^Ms, ag^Ms)<\e \} \right| <\infty .
\end{equation}
\end{defn}

Obviously this condition holds for any $g\in G$ if the action of $G$ is proper and for every loxodromic $g\in G$ if $G$ acts on $S$ acylindrically.

\begin{thm}[{\cite[Theorem 1.2]{Osi16a}}]\label{main}
For any group $G$, the following conditions are equivalent.
\begin{enumerate}
\item[(a)] $G$ is acylindrically hyperbolic.
\item[(b)] $G$ is not virtually cyclic and admits an action on a hyperbolic space such that at least one element of $G$ is loxodromic and satisfies the WPD condition.
\item[(c)] There exists a generating set $X$ of $G$ such that the corresponding Cayley graph $\Gamma (G,X)$ is hyperbolic,  $|\partial \Gamma (G,X)|> 2$, and the natural action of $G$ on $\Gamma (G,X)$ is acylindrical.
\end{enumerate}
\end{thm}

Part (c) of this theorem is especially useful for studying properties of acylindrically hyperbolic groups since it allows to pass from a (possibly non-cobounded) action of $G$ on a general hyperbolic space to the more familiar action on the Cayley graph. In addition, one can ensure that $\Gamma (G,X)$ is quasi-isometric to a tree \cite{Bal}.

\paragraph{2.3. Examples.}
Obviously every geometric action is acylindrical. In particular, this applies to the action of any finitely generated group on its Cayley graph with respect to a finite generating set. Thus every hyperbolic group is virtually cyclic or acylindrically hyperbolic. More generally, non-virtually-cyclic relatively hyperbolic groups with proper peripheral subgroups are acylindrically hyperbolic. In the latter case the action on the relative Cayley graph is non-elementary and acylindrical, see \cite{Osi16a}. Below we discuss some less obvious examples.

(a) \emph{Mapping class groups.} The mapping class group $MCG(\Sigma_{g,p})$ of a closed surface of genus $g$ with $p$ punctures is acylindrically hyperbolic unless $g=0$ and $p\le 3$ (in these exceptional cases, $MCG(\Sigma_{g,p})$ is finite). For $(g,p)\in \{ (0,4), (1,0), (1,1)\} $ this follows from the fact that $MCG(\Sigma_{g,p})$ is non-elementary hyperbolic. For all other values of $(g,p)$ this follows from hyperbolicity of the curve complex $\mathcal C(\Sigma_{g,p})$ of $\Sigma_{g,p}$ first proved by Mazur and Minsky \cite{MM} and acylindricity of the action of $MCG(\Sigma_{g,p})$ on $\mathcal C(\Sigma_{g,p})$, which is due to Bowditch \cite{Bow}.

(b) $Out(F_n)$. Let $n\ge 2$ and let $F_n$ be the free group of rank $n$. Bestvina and Feighn \cite{BFe} proved that for every fully irreducible automorphism $f\in Out(F_n)$ there exists a hyperbolic graph such that $Out(F_n)$ acts on it and the action of $f$ satisfies the weak proper discontinuity condition. Thus $Out(F_n)$ is acylindrically hyperbolic by Theorem \ref{main}.

(c) \emph{Groups acting on $CAT(0)$ spaces.} Sisto \cite{Sis} showed that if a group $G$ acts properly on a proper $CAT(0)$ space and contains a rank one element, then $G$ is either virtually cyclic or acylindrically hyperbolic. Together with the work of Caprace--Sageev \cite{CS}, this implies the following alternative for right angled Artin groups: every right angled Artin group is either cyclic, decomposes as a direct product of two non-trivial groups, or acylindrically hyperbolic. An alternative proof of the later result can be found in \cite{KK}, where Kim and Koberda construct explicitly acylindrical actions of right angled Artin groups. A similar theorem holds for graph products of groups and, even more generally, subgroups of graph products \cite{MO15}. For a survey of examples of acylindrically hyperbolic groups arising from actions on $CAT(0)$ cubical complexes, see \cite{Gen}.

(d) \emph{Fundamental groups of graphs of groups.} In \cite{MO15}, Minasyan and the author prove the following.

\begin{thm}\label{tree}
Let $G$ be a group acting minimally on a simplicial tree $T$. Suppose that $G$ does not
fix any point of $\partial T$ and there exist vertices $u,v$ of $T$ such that the pointwise stabilizer of $\{ u,v\}$ is finite. Then $G$ is either virtually cyclic or acylindrically hyperbolic.
\end{thm}

If $G$ is the fundamental group of a graph of groups $\mathcal G$, then one can apply Theorem \ref{tree} to the action of $G$ on the associated Bass-Serre tree.
In this case the minimality of the action and the absence of fixed points on $\partial T$ can be recognized from the
local structure of  $\mathcal G$. We mention here two particular cases. We say that a subgroup $C$ of a group $G$ is \emph{weakly malnormal} if there exists $g\in G$ such that $|C^g \cap C|<\infty$.

\begin{cor}\label{cor:amalg-intr}
Let $G$ split as a free product of groups $A$ and $B$ with an amalgamated subgroup $C$. Suppose $A\ne C\ne B$ and $C$ is weakly malnormal in $G$. Then $G$ is either virtually cyclic or acylindrically hyperbolic.
\end{cor}

Note that the virtually cyclic case cannot be excluded from this corollary. Indeed it realizes if $C$ is finite and has index $2$ in both factors.

\begin{cor}\label{cor:HNN-intr}
Let $G$ be an HNN-extension of a group $A$ with associated subgroups $C$ and $D$. Suppose that $C\ne A\ne D$ and $C$ is weakly malnormal in $G$. Then $G$ is acylindrically hyperbolic.
\end{cor}

These results were used in \cite{MO15} to prove acylindrical hyperbolicity of a number of groups. E.g., it implies that for every field $k$, the automorphism group $Aut\, k[x,y]$ of the polynomial algebra $k[x,y]$ is acylindrically hyperbolic. Some other applications are discussed below.

(e) \emph{$3$-manifold groups.} In the same paper \cite{MO15} (see also \cite{MOerr}), Minasyan and the author proved that for every compact orientable irreducible $3$-manifold $M$, the fundamental group $\pi_1(M)$ is either virtually polycyclic, or acylindrically hyperbolic, or $M$ is Seifert fibered. In the latter case, $\pi_1(M)$ contains a normal subgroup $N\cong \mathbb Z$ such that $\pi_1(M)/N$ is acylindrically hyperbolic.

(f) \emph{Groups of deficiency at least $2$.}  In \cite{Osi16b}, the author proved that every group which admits a finite presentation with at least $2$ more generators than relations is acylindrically hyperbolic. (The original proof contained a gap which is fixed in \cite{MOerr}.) Interestingly, the proof essentially uses results about $\ell^2$-Betti numbers of groups.

(g) \emph{Miscellaneous examples.} Other examples include central quotients of Artin-Tits groups of spherical type \cite{CW} and of $FC$ type with underlying Coxeter graph of diameter at least $3$ \cite{CM}, small cancellation groups (including infinitely presented ones)  \cite{SG}, orthogonal forms of Kac--Moody groups over arbitrary fields \cite{CH}, the Cremona group (see \cite{DGO} and references therein; the main contribution towards this result is due to Cantat and Lamy \cite{CL}), and non-elementary convergence groups \cite{Sun}.

\paragraph{2.4. Some algebraic and analytic properties.}
Our next goal is to survey some algebraic and analytic properties of acylindrically hyperbolic groups.

(a) \emph{Finite radical.} Every acylindrically hyperbolic group $G$ contains a unique maximal finite normal subgroup denoted $K(G)$ and called the \emph{finite radical} of $G$ \cite{DGO}. It also coincides with the amenable radical of $G$. In particular, $G$ has no infinite amenable normal subgroups.

(b) \emph{SQ-univerality.} Recall that a group $G$ is {\it SQ-universal} if  every countable group can be
embedded into a quotient of $G$. Informally, this property can be considered as an indication of algebraic ``largeness" of $G$. In \cite{DGO}, Dahmani, Guirardel, and the author proved the following result by using group theoretic Dehn filling (we refer to \cite{DGO} for the survey of the previous work in this direction).

\begin{thm}\label{SQ}
Every acylindrically hyperbolic group is SQ-universal.
\end{thm}

One consequence of this, also obtained in \cite{DGO}, is that every subgroup of the mapping class group $MCG(\Sigma)$ of a punctured closed surface $\Sigma$ is either  virtually abelian or SQ-universal. It is easy to show using cardinality arguments that every finitely generated SQ-universal group has uncountably many non-isomorphic quotients. This observation allows one to reprove various (well-known) non-embedding theorems for higher rank lattices in mapping class groups since these lattices have countably many normal subgroups by the Margulis normal subgroup theorem. For instance, we immediately obtain that  every homomorphism from an irreducible lattice in a connected semisimple Lie group of $\mathbb R$-rank at least $2$ with finite center to $MCG(\Sigma)$ has finite image (compare to the main result of \cite{Farb_Masur}).

(c) \emph{Mixed identities.} A group $G$ satisfies a \emph{mixed identity} $w=1$ for some $w\in G\ast F_n$, where $F_n$ denotes the free group of rank $n$, if every homomorphism $G\ast F_n\to G$ that is identical on $G$ sends $w$ to $1$. A mixed identity $w=1$ is non-trivial if $w\ne 1$ as an element of $G\ast F_n$.  We say that $G$ is \emph{mixed identity free} (or \emph{MIF} for brevity) if it does not satisfy any non-trivial mixed identity.

The property of being MIF is much stronger than being identity free and imposes strong restrictions on the algebraic structure of $G$. For example, if $G$ has a non-trivial center, then it satisfies the non-trivial mixed identity $[a,x]=1$, where  $a\in Z(G)\setminus\{ 1\}$. Similarly, it is easy to show (see \cite{HO16}) that a MIF group has no finite normal subgroups, is directly indecomposable, has infinite girth, etc. By constructing highly transitive permutation representations of acylindrically hyperbolic groups, Hull and the author proved that every acylindrically hyperbolic group with trivial finite radical is MIF \cite{HO16}.

(d) \emph{Quasi-cocycles and bounded cohomology.} The following theorem was proved in several papers under various assumptions (see \cite{BBF,BF,Ham,HO13} and references therein), which later turned out to be equivalent to acylindrical hyperbolicity.

\begin{thm}\label{Qcyc}
Suppose that a group $G$ is acylindrically hyperbolic. Let $V=\mathbb R$ or $V=\ell^p(G)$ for some $p\in [1, +\infty)$. Then the kernel of the natural map $H^2_b(G, V) \to H^2(G, V)$ is infinite dimensional. In particular, ${\rm dim\,} H^2_b(G, V)=\infty $.
\end{thm}

This result opens the door for Monod-Shalom rigidity theory for group actions on spaces with measure \cite{MS}. It also implies that acylindrically hyperbolic groups are not boundedly generated, i.e., are not products of finitely many cyclic subgroups.

(e) \emph{Stability properties.} It is not difficult to show that the class of acylindrically hyperbolic groups is stable under taking extensions with finite kernel and quotients modulo finite normal subgroups. It is also stable under taking finite index subgroups and, more generally, $s$-normal subgroups \cite{Osi16a}. Recall that a subgroup $N$ of a group $G$ is $s$-normal if $g^{-1}Ng\cap N$ is infinite for all $g\in G$.

On the other hand, it is not known if acylindrical hyperbolicity is stable under finite extensions (see \cite{MOerr}). More generally, we propose the following.

\begin{q}
\begin{enumerate}
\item[(a)] Is acylindrical hyperbolicity of finitely generated groups a quasi-isometry invariant?

\item[(b)] Is acylindrical hyperbolicity a measure equivalence invariant?
\end{enumerate}
\end{q}

The last question is partially motivated by the fact that the property $H^2_b(G, \ell^2(G))\ne 0$ enjoyed by all acylyndrically hyperbolic groups by Theorem \ref{Qcyc} is a measure equivalence invariant. For details we refer to \cite{MS}.

%%%%%%%%%%%%%%%%%%%%%%%%%%%%%%%%%%%%%%%%%%%%%%%%%%%%%%%%%%%%%%%%%%%%%%%

\section{Hyperbolically embedded subgroups}\label{3}

%%%%%%%%%%%%%%%%%%%%%%%%%%%%%%%%%%%%%%%%%%%%%%%%%%%%%%%%%%%%%%%%%%%%%%

\paragraph{3.1. Definition and basic examples.}
Hyperbolically embedded collections of subgroups were introduced in \cite{DGO} as generalizations of peripheral subgroups of relatively hyperbolic groups. To simplify our exposition we restrict here to the case of a single subgroup; the general case only differs by notation.

Let $G$ be a group, $H$ a subgroup of $G$. Suppose that $X$ is a relative generating set of $G$ with respect to $H$, i.e.,  $G=\langle X\cup H\rangle $. We denote by $\G $ the Cayley graph of $G$ whose edges are labeled by letters from the alphabet $X\sqcup H$. That is, two vertices $f,g\in G$ are connected by an edge going from $f$ to $g$ and labeled by $a\in X\sqcup H$ iff $fa=g$ in $G$. Disjointness of the union in this definition means that if a letter $h\in H$ and a letter $x\in X$ represent the same element $a\in G$, then for every $g\in G$, the Cayley graph $\G $ will have two edges connecting $g$ and $ga$: one labelled by $h$ and the other labelled by $x$.

We naturally think of the Cayley graph $\Gamma_H=\Gamma (H,H)$ of $H$ with respect to the generating set $H$ as a (complete) subgraph of $\G $.

\begin{defn}\label{he-def}
Let $G$ be a group, $H\le G$, and $X$ a (possibly infinite) subset of $G$. We say that $H$ is \emph{hyperbolically embedded in $G$ with respect to $X$} (we write $H \h (G,X)$) if $G=\langle X\cup H\rangle $ and the following conditions hold.
\begin{enumerate}
\item[(a)] The Cayley graph $\G $ is hyperbolic.
\item[(b)] For every $n\in \mathbb N$, there are only finitely many elements $h\in H$ such that the vertices $h$ and $1$ can be connected in $\G$ by a path of length at most $n$ that avoids edges of $\Gamma _H$.
\end{enumerate}
Further we say  that $H$ is hyperbolically embedded in $G$ and write $H\h G$ if $H\h (G,X)$ for some $X\subseteq G$.
\end{defn}

Note that for any group $G$ we have $G\h G$.  Indeed we can take $X=\emptyset $ in this case. Further, if $H$ is a finite subgroup of a group $G$, then $H\h G$. Indeed $H\h (G,X)$ for $X=G$. These cases are referred to as {\it degenerate}. We consider two additional examples borrowed from \cite{DGO}.

\begin{figure}
  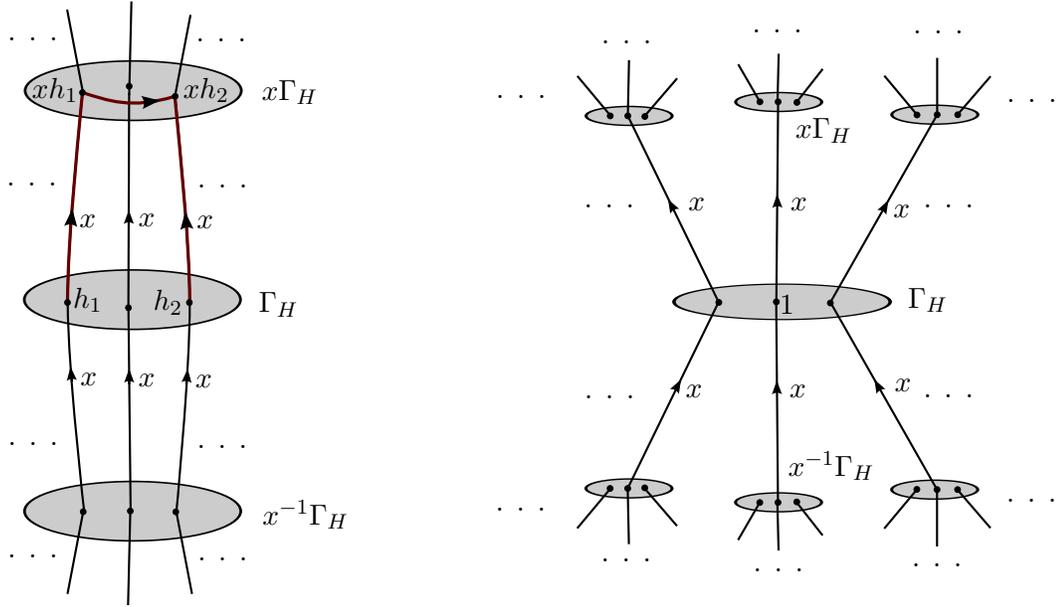
  \caption{Cayley graphs $\Gamma(G, X\sqcup H)$ for $G=H\times \mathbb Z$ and $G=H\ast \mathbb Z$.}\label{fig0}
\end{figure}

\begin{ex}\label{bex}
\begin{enumerate}
\item[(a)]
Let $G=H\times \mathbb Z$, $X=\{ x\} $, where $x$ is a generator of $\mathbb Z$. Then $\Gamma (G, X\sqcup H)$ is quasi-isometric to a line and hence it is hyperbolic. However, every two elements $h_1, h_2\in H$ can be connected by a path of length at most $3$ in $\G$ that avoids edges of $\Gamma _H$ (see Fig. \ref{fig0}). Thus $H\not\h (G,X)$ whenever $H$ is infinite.

\item[(b)]  Let $G=H\ast \mathbb Z$, $X=\{ x\} $, where $x$ is a generator of $\mathbb Z$. In this case $\Gamma (G, X\sqcup H)$ is quasi-isometric to a tree and no path connecting $h_1, h_2\in H$ and avoiding edges of $\Gamma_H$ exists unless $h_1=h_2$. Thus $H\h (G,X)$.
\end{enumerate}
\end{ex}

It is worth noting that a version of the argument from Example \ref{bex} (a) can be used to show that every hyperbolically embedded subgroup $H\h G$ is \emph{almost malnormal}, i.e., satisfies $|g^{-1}Hg\cap H|<\infty $ for all $g\in G$.

The following result is obtained in \cite{DGO} and can be regarded as a definition of relatively hyperbolic groups.

\begin{thm}\label{herh}
Let $G$ be a group, $H$ a subgroup of $G$. Then $G$ is hyperbolic relative to $H$ if and only if $H\h (G,X)$ for some finite subset $X\subseteq G$.
\end{thm}

\paragraph{3.2. Hyperbolically embedded subgroups in acylindrically hyperbolic groups.}
It turns out that acylindrical hyperbolicity of a group can be characterized by the existence of hyperbolically embedded subgroups. More precisely, the following is proved in \cite{Osi16a}.

\begin{thm}
A group $G$ is acylindrically hyperbolic if and only if it contains non-degenerate hyperbolically embedded subgroups.
\end{thm}

Moreover, in every acylindrically hyperbolic group one can find hyperbolically embedded subgroups of certain special types. We mention two results of this sort proved in \cite{DGO}. The first one plays an  important role in applications of group theoretic Dehn filling and small cancellation theory discussed below.

\begin{thm}\label{Eg}
Let $G$ be a group acting on a hyperbolic space and let $g\in G$ be a loxodromic WPD element. Then $g$ is contained in a unique maximal virtually cyclic subgroup $E(g)$ of $G$ and $E(g)\h G$.
\end{thm}

Recall that $K(G)$ denotes the final radical of an acylindrically hyperbolic group $G$ and $F_n$ denotes the free group of rank $n$.

\begin{thm}\label{FnKG}
Let $G$ be an acylindrically hyperbolic group. Then for every $n\in \mathbb N$, there exists a subgroup $H\h G$ isomorphic to $F_n\times K(G)$.
\end{thm}

The latter theorem is especially useful in conjunction with various ``extension" results proved in \cite{AHO, FPS, HO13}. Roughly speaking, these results claim that various things (e.g., group actions on metric spaces or quasi-cocycles) can be ``extended" from a hyperbolically embedded subgroup to the whole group.

%%%%%%%%%%%%%%%%%%%%%%%%%%%%%%%%%%%%%%%%%%%%%%%%%%%%%%%%%%%%%%%%%%%%%%%

\section{Group theoretic Dehn filling}\label{4}

%%%%%%%%%%%%%%%%%%%%%%%%%%%%%%%%%%%%%%%%%%%%%%%%%%%%%%%%%%%%%%%%%%%%%%

\paragraph{4.1. Dehn surgery in $3$-manifolds.}
Dehn surgery on a 3-dimensional manifold consists of cutting of a solid torus from the manifold, which may be thought of as ``drilling" along an embedded knot, and then gluing it back in a different way. The study of these ``elementary transformations" of $3$-manifolds is partially motivated by the  Lickorish-Wallace theorem, which states that every closed orientable connected 3-manifold can be obtained by performing finitely many surgeries on the $3$-dimensional sphere.

The second part of the surgery, called {\it Dehn filling}, can be formalized as follows. Let $M$ be a compact orientable 3-manifold with toric boundary. Topologically distinct ways to attach a solid torus to $\partial M$ are parameterized by free homotopy classes of unoriented essential simple closed curves in $\partial M$, called {\it slopes}. For a slope $s$, the corresponding   Dehn filling  $M(s )$ of $M$ is the manifold
obtained from $M$ by attaching a solid torus $\mathbb D^2\times \mathbb S^1$ to $\partial M$ so that the meridian
$\partial \mathbb D^2$ goes to a simple closed curve of the slope $s$.

The following fundamental theorem is due to Thurston \cite[Theorem 1.6]{Th}.

\begin{thm}[Thurston's hyperbolic Dehn surgery theorem]
Let $M$ be a compact orientable 3-manifold with toric boundary. Suppose that $M\setminus\partial M$ admits a complete finite volume hyperbolic structure. Then $M(s)$ is hyperbolic for all but finitely many slopes $s$.
\end{thm}

\paragraph{4.2. Filling in hyperbolically embedded subgroups.}
Dehn filling can be generalized in the context of abstract group theory as follows. Let $G$ be a group and let $H$ be a subgroup of $G$. One can think of $G$ and $H$ as the analogues of $\pi_1(M)$ and $\pi _1(\partial M)$, respectively. Associated to any $\sigma \in H$, is the quotient group $G/\ll s\rr $, where $\ll s\rr$ denotes the normal closure of $s$ in $G$.

If $G=\pi_1(M)$ and $H=\pi_1(\partial M)\cong \mathbb Z\oplus \mathbb Z$, where $M$ is as in Thurston's theorem, then $H$ is indeed a subgroup of $G$ and for every slope $s$, which we think of as an element of $H$, we have
\begin{equation}\label{pi1Ms}
\pi_1(M(s))=\pi_1(M)/\ll s\rr
\end{equation}
by the Seifert-van Kampen theorem. Thus $G/\ll s\rr $ is the algebraic counterpart of the filling $M(s)$.

It turns out that the analogue of Thurston's theorem holds if we start with a pair $H\le G$ such that $H$ is hyperbolically embedded in $G$. The vocabulary translating geometric terms to algebraic ones can be summarized as follows.

\medskip
\begin{center}
\begin{tabular}{|c|c|}
        \hline
       &\\
        \textbf{3-MANIFOLDS} & \textbf{GROUPS} \\
         &\\
        \hline
        \vspace{-3mm} &\\ \vspace{-3mm}
        \begin{minipage}{5.5cm}\begin{center}a compact orientable \\ 3-manifold $M$ \end{center}\end{minipage}& a group $G$ \\&\\\hline \vspace{-3mm} &\\
        \vspace{-3mm}  $\partial M$ &  $H\le G$ \\&\\\hline \vspace{-3mm} &\\
        \vspace{-3mm} \begin{minipage}{5.5cm}\begin{center}$M\setminus \partial M$ admits a finite volume \\ hyperbolic structure \end{center}\end{minipage}  & $H$ is hyperbolically embedded in $G$ \\&\\\hline \vspace{-3mm} &\\
        \vspace{-3mm} a slope $s$ & an element $h\in H$ \\&\\\hline \vspace{-3mm} &\\
        \vspace{-3mm}
        $M(s )$ & $G/\ll h\rr$
        \\&\\\hline
 \end{tabular}
\end{center}

\medskip

In these settings, the analogue of Thurston's theorem was proved in \cite{DGO}. Note that instead of considering single elements of $H$, we allow normal subgroups generated by arbitrary sets of elements. A number of additional properties can be added to the main statements (a)--(c); we mention just one of them, which is necessary for the applications considered in the next section.

\begin{thm}\label{CEP}
Let $G$ be a group, $H$ a subgroup of $G$. Suppose that $H\h (G,X)$ for some $X\subseteq G$. Then there exists a finite subset $\mathcal F$ of nontrivial elements of $H$ such that for every subgroup $N\lhd H$ that does not contain elements of $\mathcal F$, the following hold.
\begin{enumerate}
\item[(a)] If $G$ is acylindrically hyperbolic, then so is $G/\ll N\rr $, where $\ll N\rr $ denotes the normal closure of $N$ in $G$.

\item[(b)] The natural map from $H /N$ to $G/\ll N\rr $ is injective (equivalently, $H\cap \ll N\rr =N$).

\item[(c)] $H/N\h (G/\ll N\rr , \overline{X})$, where $\overline{X}$ is the natural image of $X$ in $G/\ll N\rr$.

\item[(d)] $\ll N\rr $ is the free product of conjugates of $N$ in $G$ and every element of $\ll N\rr $ is either conjugate to an element of $N$ or acts loxodromically on $\Gamma (G, X\sqcup H)$.
\end{enumerate}
\end{thm}

Note that if $H\h G$ is non-degenerate, then $G$ is always acylindrically hyperbolic. However the theorem holds (trivially) for degenerate hyperbolically embedded subgroups as well.

Combining this theorem with Theorem \ref{herh} and some basic properties of relatively hyperbolic groups, we obtain the following result, which was first proved by the author in \cite{Osi07}. It was also independently proved by Groves and Manning \cite{GM} under the additional assumptions that the group $G$ is torsion free and finitely generated.

\begin{cor} \label{CEPrh} Suppose that a group $G$ is hyperbolic relative to a subgroup $H\ne G$. Then for any subgroup $N\lhd H$ avoiding a fixed finite set of nontrivial elements, the natural map from $H/N$ to $G/\ll N\rr $ is injective and $G/\ll N\rr $ is hyperbolic relative to $H/N$. In particular, if $H/N$ is hyperbolic, then so is $G/\ll N\rr $; if, in addition, $G$ is non-virtually-cyclic, then so is $G/\ll N\rr $.
\end{cor}

Under the assumptions of Thurston's theorem, we have $H=\pi_1(\partial M)=\mathbb Z\oplus\mathbb Z$. Slopes in  $\partial M$ correspond to non-trivial primitive elements $s\in H$; for every such $s$, we have $H/\langle s\rangle \cong \mathbb Z$. Applying Corollary \ref{CEPrh} to $N=\langle s\rangle\lhd H$, we obtain that $G/\ll N\rr $ is not virtually cyclic and hyperbolic. Modulo the geometrization conjecture this algebraic statement is equivalent to hyperbolicity of $M(s)$. Thus parts (a)--(c) of Theorem \ref{CEP} indeed provide a group theoretic generalization of Thurston's theorem.

\paragraph{4.3. Applications.}
It is not feasible to discuss all applications of group theoretic Dehn surgery in a short survey. Here we list some of the results which make use of Theorem \ref{CEP} or its relatively hyperbolic analogue, Corollary \ref{CEPrh}, and provide references for further reading. We then pick one application and discuss it in more detail.

(a)\hspace{2mm} \emph{The virtual Haken conjecture.} Group theoretic Dehn filling in relatively hyperbolic groups, along with Wise's machinery of virtually special groups, was used in Agol's proof of the virtual Haken conjecture \cite{A}. Additional results on Dehn filling necessary for the proof were obtained by Agol, Groves, and Manning in the appendix to \cite{A}. One piece of Wise's work used in \cite{A} is the malnormal special quotient theorem; Agol, Groves, and Manning also found an alternative proof of this result based on Dehn filling technique \cite{AGM}.

(b)\emph{The isomorphism problem for relatively hyperbolic groups.} In \cite{DG,DT}, Dahmani, Guirardel, and Touikan, used Dehn filling to solve the isomorphism problem for relatively hyperbolic groups with residually finite parabolic subgroups under certain additional assumptions. The main idea is to apply (an elaborated version of) Corollary \ref{CEPrh} and some other results from \cite{DGO} to finite index normal subgroups in parabolic groups. This yields an approximation of relatively hyperbolic groups by hyperbolic ones, which in turn allows the authors make use of the solution of the isomorphism problem for hyperbolic groups obtained in \cite{DG11}.

(c) \emph{Residual finiteness of outer automorphism groups.} In \cite{MO10}, Dehn filling in relatively hyperbolic groups was used by Minasyan and the author to prove that $Out(G)$ is residually finite for every residually finite group $G$ with infinitely many ends; in general, this result fails for one ended groups. Results of \cite{MO10} were recently generalized to acylindrically hyperbolic groups by Antolin, Minasyan, and Sisto. In particular, they proved residual finiteness of mapping class groups of certain Haken $3$-manifolds. Acylindrical hyperbolicity of $3$-manifold groups plays a crucial role in the proof.

(d) \emph{Primeness of von Neumann algebras.} Chifan, Kida, and Pant \cite{CKP} used Dehn filling to prove primeness of von Neumann algebras of certain relatively hyperbolic groups.

(e) \emph{Farell-Jones conjecture for relatively hyperbolic groups.} Bartels \cite{B} proved that the class of groups satisfying the Farell-Jones conjecture is stable under relative hyperbolicity. In the particular case when peripheral subgroups are residually finite, an alternative proof based on Dehn filling was found by Antolin, Coulon, and Gandini \cite{ACG}.

(f) \emph{SQ-universality of acylindrically hyperbolic groups.} One simple application of Theorem \ref{CEP} is the proof of Theorem \ref{SQ}. It follows easily from $SQ$-universality of free groups of rank $2$, Theorem \ref{FnKG}, and part (b) of Theorem \ref{CEP}. For details, see \cite{DGO}.

\paragraph{4.4. Purely pseudo-Anosov subgroups of mapping class groups.}
We illustrate Theorem \ref{CEP} by considering an application to mapping class groups. Recall that a subgroup of a mapping class group is called \emph{purely pseudo-Anosov}, if all its non-trivial elements are pseudo-Anosov. The following question is Problem 2.12(A) in  Kirby's list: \emph{Does the mapping class group of any  closed orientable surface of genus $g\ge 1$  contain a non-trivial purely pseudo-Anosov normal subgroup?} It was asked in the early 1980s and is often attributed to Penner,  Long, and  McCarthy. It is also recorded by  Ivanov \cite[Problems 3]{Iv}, and  Farb refers to it in \cite{F_book} as a ``well known open question".

The abundance of finitely generated non-normal purely pseudo-Anosov free subgroups of mapping class groups is well known, and follows from an easy ping-pong argument. However, this method does not allow one to construct normal subgroups, which are usually infinitely generated. For a surface of genus $2$ the question was answered by  Whittlesey \cite{Whi} who proposed an example based on Brunnian braids. Unfortunately the methods of \cite{Whi} do not generalize even to closed surfaces of higher genus.

Another question was probably first asked by Ivanov (see \cite[Problem 11]{Iv}): \emph{Is the normal closure of a certain nontrivial power of a pseudo-Anosov element of $MCG(S_g)$ free?} Farb also recorded this question in \cite[Problem 2.9]{F_book}, and qualified it as a ``basic test question" for understanding normal subgroups of mapping class groups.

We answer both questions positively. In fact, our approach works in more general settings.

\begin{thm}[Theorem 2.30, \cite{DGO}]\label{wpd-free}
Let $G$ be a group acting on a hyperbolic space $S$, $g\in G$ a WPD loxodromic element.  Then there exists $n\in \mathbb N$ such that the normal closure $\ll g^n\rr$ in $G$ is free and purely loxodromic, i.e., every nontrivial element of $\ll g^n\rr$ acts loxodromically on $S$.
\end{thm}

This result can be viewed as a generalization of Delzant's theorem \cite{Del} stating that for a hyperbolic group $G$ and every element of infinite order $g\in G$, there exists $n\in \mathbb N$ such that $\ll g^n\rr $ is free (see also \cite{Chay} for a clarification of certain aspects of Delzant's proof).

The idea of the proof is the following. By Theorem \ref{Eg}, $g$ is contained in the maximal virtually cyclic subgroup $E(g)$ which is hyperbolically embedded in $G$. Since $\langle g\rangle $ has finite index in $E(g)$, we have $\langle g^n\rangle \lhd E(g)$. Passing to a multiple of $n$ if necessary, we can ensure that $\langle g^n\rangle$ avoids any finite collection of non-trivial elements. Thus we can apply Theorem \ref{CEP} to $H=E(g)$ and $N=\langle g^n\rangle$. Since $\langle g^n\rangle \cong \mathbb Z$, part (d) of the theorem implies that $\ll g^n\rr$ is free. That $\ll g^n\rr$ is purely loxodromic also follows from part (d) and some additional arguments relating $\G$ to $S$.

Applying  Theorem \ref{wpd-free} to mapping class groups acting on the curve complexes, we obtain the following.

\begin{cor}
Let $\Sigma$ be a possibly punctured closed orientable surface. Then for any pseudo-Anosov element  $a\in MCG(\Sigma)$, there exists $n\in \mathbb N$ such that the normal closure of $a^{n}$ is free and purely pseudo-Anosov.
\end{cor}

\section{Small cancellation theory and its applications}

\paragraph{5.1. Generalizing classical small cancellation.}
The classical small cancellation theory deals with presentations
$$
F(X)/\ll \mathcal R\rr=\langle X\mid \mathcal R\rangle,
$$
where $F(X)$ is the free group with basis $X$, and common subwords of distinct relators are ``small" in a certain precise sense. This property allows one to control cancellation in products of conjugates of relators (and their inverses); in turn, this leads to nice structural results for the normal closure $\ll \mathcal R\rr$ and the group $F(X)/\ll \mathcal R\rr$.

More generally, one can replace the free group $F(X)$ with a group $G_0$ enjoying some hyperbolic properties and add new relations to a presentation of $G_0$. If these new relations satisfy a suitable version of small cancellation, many results of the classical small cancellation theory can be reproved in these settings. On the other hand, the small cancellation assumptions are usually general enough to allow one to create interesting relations between elements.

The idea of generalizing classical small cancellation to groups acting on hyperbolic spaces is due to Gromov \cite{Gro}, although some underlying ideas go back to the work of Olshanskii \cite{Ols82,Ols80}.  In the case of hyperbolic groups, it was formalized by Delzant \cite{Del}, Olshanskii \cite{Ols93}, and others. Olshanskii's approach was generalized to relatively hyperbolic groups by the author in \cite{Osi10} and further generalized to acylindrically hyperbolic groups by Hull \cite{Hull}. These generalizations employ isoperimetric characterizations of relatively hyperbolic groups and hyperbolically embedded subgroups \cite{DGO,Osi06a} and follow closely the classical theory.
Yet another approach is based on Gromov's \emph{rotating families} (see \cite{DGO} and references therein.)

Unfortunately, the ideas involved in this work are too technical for a short survey paper and we do not discuss them here. Instead we discuss one applications of small cancellation theory in relatively hyperbolic groups to proving embedding theorems and studying conjugacy growth of groups \cite{Osi10,HO11}.

\paragraph{5.2. Embedding theorems and conjugacy growth of groups.}
In 1949, Higman, B.H. Neumann, and H. Neumann
proved that any countable group $G$ can be embedded into a
countable group $B$ such that every two elements of the same order
are conjugate in $B$ \cite{HNN}. We notice that the group $B$ in \cite{HNN}
is constructed as a union of infinite number of subsequent
HNN--extensions and thus $B$ is never finitely generated. In \cite{Osi10}, the author used small cancellation theory in relatively hyperbolic groups to prove the following stronger result. For a group $G$, we
denote by $\pi (G)$ the set of finite orders of elements of
$G$.

\begin{thm}\label{Conj}
Any countable group $G$ can be embedded into a  finitely generated
group $C$ such that any two elements of the same order are
conjugate in $C$ and $\pi (G)=\pi (C)$.
\end{thm}

We explain the idea of the proof in the particular case when $C$ is torsion free. Let $G_0=C\ast F(x,y)$, where $F(x,y)$ is the free group with basis $\{ x, y\}$. Given any non-trivial lement $g\in G_0$, one first considers the HNN-extension $$H=\langle G_0, t\mid t^{-1}gt=x\rangle .$$ Obviously $x$ and $g$ are conjugate in $H$. Then imposing an additional relation $t=w(x,y)$, where $w(x,y)$ is a suitable small cancellation word in the alphabet $\{ x,y\}$, one ensures that this conjugation happens in a certain quotient group $G_1$ of $G_0$. Small cancellation theory is then used to show that the restriction of the natural homomorphism $G_0\to G_1$ to $C$ is injective and the image of $F(x,y)$ in $G_1$ is still ``large enough". Here ``large enough" means that the image of $F(x,y)$ in $G_1$ is non-elementary with respect to some acylindrical action of $G_1$ on a hyperbolic space. This allows us to iterate the process. Repeating it for all non-trivial elements we obtain a group with $2$ conjugacy classes which is generated by $2$ elements (the images of $x$ and $y$) and contains $C$.

Applying Theorem \ref{Conj} to the group $G=\mathbb Z$, we obtain the following.

\begin{cor}
There exists a torsion free finitely generated group with $2$ conjugacy classes.
\end{cor}

The existence of a finitely generated group with $2$ conjugacy classes other than $\mathbb Z/2\mathbb Z$ was a long standing open problem, sometimes attributed to Maltsev. It is easy to see that such groups do not exist among finite (and residually finite) groups. It is also observed in \cite{Osi10} that such a group cannot be constructed as a limit of hyperbolic groups; this justifies the use of small cancellation theory in the more general settings.

Given a group $G$ generated by a finite set $X$, the associated {\it conjugacy growth function} of $G$, denoted by $\xi_{G,X} $, is defined as follows: $\xi _{G,X}(n)$ is the number of conjugacy classes of elements that can be represented by words of length at most $n$ is the alphabet $X\cup X^{-1}$. Given $f,g\colon\mathbb N \to \mathbb N $, we write $f\sim g$ if there exists $C\in \mathbb N$ such that  $f(n)\le g(Cn)$ and $g(n)\le f(Cn)$ for all $n\in \mathbb N$. Obviously $\sim $ is an equivalence relation and $\xi _{G,X}(n)$ is independent of the choice of $X$ up to this equivalence.

The conjugacy growth function was introduced by Babenko \cite{IB} in order to study geodesic growth of Riemannian manifolds. For more details and a survey of some recent results about conjugacy growth we refer to \cite{HO11}. Based on ideas from \cite{Osi10}, Hull and the author also obtained a complete description of functions that occur as conjugacy growth functions of finitely generated groups.  It is worth noting that such a description for the usual growth function seems to be out of reach at this time.

\begin{thm}
Let $G$ be a group generated by a finite set $X$, and let $f$ denote the conjugacy growth function of $G$ with respect to $X$. Then the following conditions hold.
\begin{enumerate}
\item[(a)] $f$ is non-decreasing.
\item[(b)] There exists $a\ge 1$ such that $f(n) \le  a^n$ for every $n\in\mathbb N$.
\end{enumerate}
Conversely, suppose that a function $f\colon\mathbb N\to \mathbb N$ satisfies the above conditions (a) and (b). Then there exists an group $G$ generated by a finite set $X$ such that $\xi_{G,X}\sim f$.
\end{thm}

Of course, the non-trivial part of the theorem is the fact that every function satisfying (a) and (b) realizes as the conjugacy growth function.

Yet another result proved in \cite{HO11} is the following.

\begin{thm}
There exists a finitely generated group $G$ and a finite index subgroup $H\le G$ such that $H$ has $2$ conjugacy classes while $G$ is of exponential conjugacy growth.
\end{thm}

In particular, unlike the usual growth function, conjugacy growth of a group is not a quasi-isometry invariant.

Readers interested in other applications of small cancellation technique to groups with hyperbolically embedded subgroups are referred to \cite{Hull} and \cite{MO18}; for a slightly different approach employing rotating families see Gromov's paper \cite{Gro03}, Coulon's survey \cite{Coulon}, and references therein.


\begin{thebibliography}{99}

\bibitem{AHO}
C. Abbott, D. Hume, D. Osin, Extending group actions on metric spaces, arXiv:1703.03010.

\bibitem{A}
I. Agol, The virtual Haken conjecture. With an appendix by I. Agol, D. Groves, and J. Manning. \emph{Doc. Math.} \textbf{18} (2013), 1045-1087.

\bibitem{AGM}
I. Agol, D. Groves, J. Manning, An alternate proof of Wise's malnormal special quotient theorem,
\emph{Forum Math. Pi} \textbf{4} (2016), e1.

\bibitem{ACG}
Y. Antolin, R. Coulon, G. Gandini, Farrell-Jones via Dehn fillings, arXiv:1510.08113.


\bibitem{AMO} G. Arzhantseva, A. Minasyan, D. Osin, The SQ-universality and residual properties of relatively hyperbolic groups, {\it J. Algebra} {\bf 315} (2007), no. 1, 165-177.

\bibitem{Bal}
S. Balasubramanya, Acylindrical group actions on quasi-trees, arXiv:1602.03941.


\bibitem{B}
A. Bartels, Coarse flow spaces for relatively hyperbolic groups, \emph{Compositio Math.} \textbf{153} (2017), 745-779.

\bibitem{BBF}
M. Bestvina, K. Bromberg, K. Fujiwara, Bounded cohomology with coefficients in uniformly convex Banach spaces, \emph{Comment. Math. Helv.} \textbf{91} (2016), no. 2, 203-218.

\bibitem{BFe}
M. Bestvina, M. Feighn, A hyperbolic $Out(F_n)$-complex,\emph{ Groups Geom. Dyn.} \textbf{4} (2010), no. 1, 31-58.

\bibitem{BF}
M. Bestvina, K. Fujiwara, Bounded cohomology of subgroups of mapping class groups, \emph{Geom. Topol.} \textbf{6} (2002), 69–89.

\bibitem{Bow}
B. Bowditch, Tight geodesics in the curve complex, \emph{Invent. Math.} \textbf{171} (2008), no. 2, 281-300.

\bibitem{BH} M. Bridson, A. Haefliger, Metric spaces of non-positive curvature.
Grundlehren der Mathematischen Wissenschaften, 319. Springer-Verlag, Berlin, 1999.

\bibitem{IB} I. Babenko, Closed geodesics, asymptotic volume, and characteristics of group growth, \textit{Izv. Akad. Nauk SSSR, Ser. Mat.} \textbf{52} (1988), 675-711.

\bibitem{CW} M. Calvez, B. Wiest, Acylindrical hyperbolicity and Artin-Tits groups of spherical type, 	 arXiv:1606.07778.

\bibitem{CL}
S. Cantat, S. Lamy,
Normal subgroups in the Cremona group. With an appendix by Y. de Cornulier.
\textit{Acta Math.} \textbf{210} (2013) no. 1, 31--91.

\bibitem{CH}
P.-E. Caprace, D. Hume, Orthogonal forms of Kac--Moody groups are acylindrically hyperbolic, 	arXiv:1408.6117.

\bibitem{CS}
P.-E. Caprace, M. Sageev,
Rank rigidity for $CAT(0)$ cube complexes,
\textit{Geom. Funct. Anal.}
\textbf{21} no. 4 (2011), 851-891.

\bibitem{CM}
I. Chatterji, A. Martin
A note on the acylindrical hyperbolicity of groups acting on $CAT(0)$ cube complexes, arXiv:1610.06864


\bibitem{Chay}
V. Chaynikov, On the generators of the kernels of hyperbolic group presentations,
\emph{Algebra Discrete Math.} \textbf{11} (2011), no. 2, 18-50.

\bibitem{CKP}
I. Chifan, Y. Kida, S. Pant, Primeness results for von Neumann algebras associated with surface braid groups, \emph{Int. Math. Res. Not.} (2016), no. 16, 4807-4848.

\bibitem{Coulon}
R. Coulon, Th\'eorie de la petite simplification: une approche g\'eom\'etrique [d'apr\`es F. Dahmani, V. Guirardel, D. Osin et S. Cantat, S. Lamy] (in French), Ast\'erisque No. 380, S\'eminaire Bourbaki, Vol. 2014/2015 (2016), Exp. No. 1089, 1-33.

\bibitem{DG}
F. Dahmani, V. Guirardel, Recognizing a relatively hyperbolic group by its Dehn fillings, arXiv:1506.03233.


\bibitem{DG11}
F. Dahmani, V. Guirardel, The isomorphism problem for all hyperbolic groups, \emph{Geom. Funct. Anal.} \textbf{21} (2011), no. 2, 223–300.
\bibitem{DGO} F. Dahmani, V. Guirardel, D. Osin, Hyperbolically embedded subgroups and rotating families in groups acting on hyperbolic spaces,  \emph{Memoirs AMS}  \textbf{245} (2017), no. 1156.

\bibitem{DT}
F. Dahmani, N. Touikan,
Deciding Isomorphy using Dehn fillings, the splitting case, arXiv:1311.3937.


\bibitem{Del} T. Delzant,
Sous-groupes distingu\'es et quotients des groupes hyperboliques,
\textit{Duke Math. J. }
\textbf{83}  (1996),  no. 3, 661-682.

\bibitem{F_book}
B. Farb,
Some problems on mapping class groups and moduli space,
in \textit{ Problems on mapping class groups and related topics},
11-55, Proc. Sympos. Pure Math., 74, Amer. Math. Soc., Providence, RI.

\bibitem{Farb_Masur}
B. Farb, H. Masur,
Superrigidity and mapping class groups,
\textit{Topology}
\textbf{37} (1998), no. 6, 1169-1176.

\bibitem{FPS}
R. Frigerio, M. Pozzetti, A. Sisto, Extending higher-dimensional quasi-cocycles, \emph{J. Topol.} \textbf{8} (2015), no. 4, 1123-1155.

\bibitem{Gen}
A. Genevois, Hyperbolicities in CAT(0) cube complexes, arXiv:1709.08843.

\bibitem{GH}
E. Ghys, P. de la Harpe, Sur les groupes hyperboliques d'apr\`es Mikhael Gromov. Progress in Mathematics, \textbf{83}. Birkh\"auser Boston, Inc., Boston, MA, 1990.

\bibitem{Gro03}
M. Gromov, Random walk in random groups, \emph{Geom. Funct. Anal.} \textbf{13} (2003), no. 1, 73-146.

\bibitem{Gro}
M. Gromov, Hyperbolic groups, Essays in Group Theory, MSRI Series,
Vol.8, (S.M. Gersten, ed.), Springer, 1987, 75--263.

\bibitem{GM} D. Groves, J. Manning, Dehn filling in relatively hyperbolic groups, \emph{Israel J. Math.} \textbf{168} (2008), 317-429.

\bibitem{SG}
D. Gruber, A. Sisto,
Infinitely presented graphical small cancellation groups are acylindrically hyperbolic, arXiv:1408.4488.


\bibitem{Ham}
U. Hamenst\"adt, Bounded cohomology and isometry groups of hyperbolic spaces,
\emph{J. Eur. Math. Soc.} \textbf{10} (2008), no. 2, 315-349.

\bibitem{HNN}
G. Higman, B.H. Neumann, H. Neumann,  Embedding theorems for
groups, \emph{J. London Math. Soc.} {\bf 24} (1949), 247-254.

\bibitem{Hull}
M. Hull, Small cancellation in acylindrically hyperbolic groups, \emph{Groups Geom. Dyn.} \textbf{10} (2016), no. 4, 1077-1119.

\bibitem{HO16} M. Hull, D. Osin, Transitivity degrees of countable groups and acylindrical hyperbolicity, \emph{Israel J. Math.} \textbf{216} (2016), no. 1, 307–353.

\bibitem{HO13}
M. Hull, D. Osin, Induced quasi-cocycles on groups with hyperbolically embedded subgroups, \emph{Alg. \& Geom. Topol.}, \textbf{13} (2013) 2635-2665.

\bibitem{HO11}
M. Hull, D. Osin, Conjugacy growth of finitely generated groups, \emph{Adv. Math.}  \textbf{235} (2013), 361-389; Corrigendum to ``Conjugacy growth of finitely generated groups", \emph{Adv. Math.} \textbf{294} (2016), 857-859.

\bibitem{Iv}
N. Ivanov,
Fifteen problems about the mapping class groups,
in \textit{Problems on mapping class groups and related topics},
71-80, Proc. Sympos. Pure Math., 74, Amer. Math. Soc., Providence, RI, 2006.


\bibitem{KK}
S.-H. Kim, T. Koberda, The geometry of the curve graph of a right-angled Artin group,
\emph{Internat. J. Algebra Comput.} \textbf{24} (2014), no. 2, 121-169.

\bibitem{MM}
H. Masur, Y. Minsky, Geometry of the complex of curves. I. Hyperbolicity.
\emph{Invent. Math.} \textbf{138} (1999), no. 1, 103-149.

\bibitem{MO18}
A. Minasyan, D. Osin, Acylindrically hyperbolic groups with exotic properties,  arXiv:1804.08767.

\bibitem{MOerr}
A. Minasyan, D. Osin, Erratum to the paper "Acylindrical hyperbolicity of groups acting on trees", arXiv:1711.09486.


\bibitem{MO15} A. Minasyan, D. Osin, Acylindrical hyperbolicity of groups acting on trees, \emph{Math. Ann.} \textbf{362} (2015), no. 3-4, 1055-1105.

\bibitem{MO10} A. Minasyan, D. Osin, Normal automorphisms of relatively hyperbolic groups, \emph{Trans. AMS}, {\bf 362} (2010), no. 11, 6079-6103.

\bibitem{MS}
N. Monod, Y. Shalom,
Orbit equivalence rigidity and bounded cohomology,\emph{ Ann. of Math.} \textbf{164} (2006), no. 3, 825-878.

\bibitem{Ols93}
A. Olshanskii,  On residualing homomorphisms and
$G$--subgroups of hyperbolic groups, \textit{Internat. J.
Algebra Comput.} {\bf 3}
(1993), 4, 365-409.

\bibitem{Ols82}
A. Olshanskii, The Novikov-Adyan theorem (in Russian), \emph{Mat. Sb.} \textbf{118}(160) (1982), no. 2, 203-235.

\bibitem{Ols80}
A. Olshanskii, An infinite group with subgroups of prime orders (in Russian), \emph{Izv. Akad. Nauk SSSR Ser. Mat.} \textbf{44} (1980), no. 2, 309-321.

\bibitem{Osi16a} D. Osin, Acylindrically hyperbolic groups, \emph{Trans. Amer. Math. Soc.} \textbf{368} (2016), no. 2, 851-888.

\bibitem{Osi16b} D. Osin, On acylindrical hyperbolicity of groups with positive first $\ell^2$-Betti number,  \emph{Bull. Lond. Math. Soc.} \textbf{47} (2015), no. 5, 725-730.

\bibitem{Osi10} D. Osin, Small cancellations over relatively hyperbolic groups and embedding theorems, \emph{Ann. Math.} {\bf 172} (2010), no. 1, 1-39.

\bibitem{Osi07}  D. Osin, Peripheral fillings of relatively hyperbolic groups, \textit{Invent. Math.} {\bf 167} (2007), no. 2, 295-326.

\bibitem{Osi06a} D. Osin, Relatively hyperbolic groups: Intrinsic
geometry, algebraic properties, and algorithmic problems, {\it
Memoirs AMS} {\bf 179} (2006), no. 843, vi+100 pages.


\bibitem{Sel}
Z. Sela, Acylindrical accessibility for groups,
\emph{Invent. Math.} \textbf{129} (1997), no. 3, 527-565.

\bibitem{Sis}
A. Sisto, Contracting elements and random walks, arXiv:1112.2666.

\bibitem{Sun}
B. Sun, A dynamical characterization of acylindrically hyperbolic groups, arXiv:1707.04587.

\bibitem{Th}
W.P. Thurston,
Three-dimensional manifolds, Kleinian groups and hyperbolic geometry,
{\it Bull. Amer. Math. Soc. }
{\bf 6} (1982), no. 3, 357-381.

\bibitem{Whi}
K. Whittlesey,
Normal all pseudo-Anosov subgroups of mapping class groups,
\textit{Geom. Topol.}
\textbf{4}  (2000), 293--307.

\end{thebibliography}
\end{document}